\documentclass[11pt,twoside,a4paper]{article}
\usepackage{amsmath,amssymb,amsthm}

\usepackage{graphicx,psfrag}
\newtheorem{Thm}{Theorem}[section]

\theoremstyle{definition}

\theoremstyle{remark}
\newtheorem{Rem}[Thm]{Remark}

\newcommand{\N}{\mathbb{N}}

\newcommand{\ep}{\varepsilon}

\newcommand{\si}{\sigma}

\renewcommand{\phi}{\varphi}

\begin{document}

\title{Quasi-metric and metric spaces}
\author{Viktor Schroeder\footnote{Supported by Swiss National
Science Foundation}}

\date{}
\maketitle

\begin{abstract} 
We give a short review of a construction of Frink to
obtain a metric space from a quasi-metric space.
By an example we illustrate the limits of the construction.

\end{abstract}

\section{Introduction}

A {\em quasi-metric space}
is a set
$Z$
with a function
$\rho:Z\times Z\to [0,\infty)$
which satisfies the conditions
\begin{itemize}
\item[(1)] $\rho(z,z')\ge 0$
for every
$z$, $z'\in Z$
and
$\rho(z,z')=0$
if and only if
$z=z'$;
\item[(2)]
$\rho(z,z')=\rho(z',z)$
for every
$z$, $z'\in Z$;
\item[(3)]
$\rho(z,z'')\le K\max\{\rho(z,z'),\rho(z',z'')\}$
for every
$z$, $z'$, $z''\in Z$
and some fixed
$K\ge 1$.

\end{itemize}

The function
$\rho$
is called in that case a
{\em quasi-metric},
or more specifically, a
$K$-quasi-metric.
The property (3) is a generalized version of
the {\em ultra-metric}
triangle inequality (the case
$K=1$).

\begin{Rem}\label{remark:powermetric}
If
$(Z,d)$
is a metric space, then
$d$
is a
$K$-quasi-metric for
$K=2$.
In general
$d^p$
is not a metric on
$Z$
for
$p>1$.
But
$d^p$
is still a
$2^p$-quasi-metric.
\end{Rem}

We are interested in the question how to obtain a metric
on $Z$. Our personal motivation comes from the study of the boundary
at infinity of a Gromov hyperbolic space, where this question arises naturally,
see e.g.
\cite[chapter 2]{BS}.
The problem was studied by Frink in the interesting paper \cite{Fr}. The motivation of
Frink was to obtain suitable condition for a topological space to be
metrizable.
Frink used a natural approach, which we call the
{\em chain approach} to obtain a metric on $Z$.
He showed that this appraoch works and gives a metric if 
the space
$(Z,\rho)$ satisfies 
the axioms (1), (2) from above and 
instead of (3) the {\em weak triangle inequality}

(3') If $\rho(z,z') \le \ep$ and $\rho(z',z'') \le \ep$
then $\rho(z,z'') \le 2\ep$.

\noindent Observe that (3') is equivalent to (3) with constant $K=2$ but the
formulation as a weak triangle inequality 
points out that the constant $K=2$ playes a special role.

In this short note we give a review of Frinks approach
and show that 
there exists 
a "natural" counterexample to the chain approach in case
that the weak triangle inequality (3') is not satisfied.

I am grateful to Urs Lang for showing me the paper of Frink and
to Tobias Strubel for discussions about quasimetric spaces.

\subsection{Quasimetrics and metrics}

Let 
$(Z,\rho)$ be a quasimetric space.
We want to obtain a metric
on $Z$. 
Since $\rho$ satisfies all axioms of a metric space except 
the triangle inequality, the following 
approach is very natural.
Define a map
$d:Z\times Z \to [0,\infty)$,
$d(z,z')=\inf\sum_i\rho(z_i,z_{i+1}),$
where the infimum is taken over all sequences
$z=z_0,\dots,z_{n+1}=z'$
in
$Z$.
By definition $d$ satisfy the
triangle inequality. 
We call this approach to the triangle
inequality the 
{\it chain approach}.
The problem with the chain approach is that
$d(z,z')$ could be $0$ for different
points $z, z'$ and axiom
(1) is not longer satisfied 
for $(Z,d)$.

This chain approach is due to
Frink, who realized
that the approach works if
the space
$(Z,\rho)$ satisfies 
the axioms (1), (2) and (3') from above. 
For the convenience of the reader we give a proof
of Frink's result.

\begin{Thm}\label{Thm:chainconstract} Let
$\rho$
be a
$K$-quasi-metric
on a set
$Z$
with
$K\le 2$.
Then, the chain construction applied to
$\rho$
yields a metric
$d$
with
$\frac{1}{2K}\rho\le d\le\rho$.
\end{Thm}

\begin{proof} Clearly,
$d$
is nonnegative, symmetric, satisfies the
triangle inequality and
$d\le\rho$.
We prove by induction over the length of sequences
$\si=\{z=z_0,\dots,z_{k+1}=z'\}$, $|\si|=k+2$,
that
\begin{equation}\label{eqn:qmetricest}
   \rho(z,z')\le\sum(\si):=
   K\left(\rho(z_0,z_1)+2\sum_1^{k-1}\rho(z_i,z_{i+1})
   +\rho(z_k,z_{k+1})\right).
\end{equation}
For
$|\si|=3$,
this follows from the triangle inequality (3) for
$\rho$.
Assume that (\ref{eqn:qmetricest}) holds true for
all sequences of length
$|\si|\le k+1$,
and suppose that
$|\si|=k+2$.

Given
$p\in\{1,\dots,k-1\}$,
we let
$\si_p'=\{z_0,\dots,z_{p+1}\}$,
$\si_p''=\{z_p,\dots,z_{k+1}\}$,
and note that
$\sum(\si)=\sum(\si_p')+\sum(\si_p'')$.

Because
$\rho(z,z')\le K\max\{\rho(z,z_p),\rho(z_p,z')\}$,
there is a maximal
$p\in\{0,\dots,k\}$
with
$\rho(z,z')\le K\rho(z_p,z')$.
Then
$\rho(z,z')\le K\rho(z,z_{p+1})$.

Assume now that
$\rho(z,z')>\sum(\si)$.
Then, in particular,
$\rho(z,z')>K\rho(z,z_1)$
and
$\rho(z,z')>K\rho(z_k,z')$.
It follows that
$p\in\{1,\dots,k-1\}$
and thus by the inductive assumption
$$\rho(z,z_{p+1})+\rho(z_p,z')\le
  \sum(\si_p')+\sum(\si_p'')=\sum(\si)<\rho(z,z').$$
On the other hand,
$$\rho(z,z')\le K\min\{\rho(z,z_{p+1}),\rho(z_p,z')\}
  \le\rho(z,z_{p+1})+\rho(z_p,z')$$
because
$K\le 2$,
a contradiction. Now, it follows from (\ref{eqn:qmetricest})
that
$\rho\le 2Kd$,
hence,
$d$
is a metric as required.
\end{proof}

%%%%%%%%%%%%%%%%%%%%%%%%%%%%%%%%%%%%%%%%%%%%%%%%%%%%%%%%%%%%%%%%%%%%%%

\section{Example}

In this section we construct for any given 
$\ep > 0$ a quasimetric space $(Z,\rho)$ 
such that the chain approach does not lead to a metric space and
the following holds:
for every triple of points $z_0,z_1,z_2$ we have
$$\rho(z_1,z_2) \le (1+\ep) [\rho(z_1,z_0)+\rho(z_0,z_2)].$$

Let therefore $a \in (0,\frac{1}{2})$ be a given constant.
Let $Z$ be the set of dyadic rational of the interval $[0,1]$.
Then $Z$ is the disjoint union of
$Z_n$ , $n\in \N$, where $Z_0=\{0,1\}$, and 
$Z_n=\{\frac{k}{2^n} : 0 < k < 2^n, \mbox{{\em k odd}}\}$
for $n\ge 1$.
If $z \in Z_n$, we say that the level of
$z$ is $n$ and write
$\ell(z)=n$.
For the following construction it is useful to see
$Z$ embedded by
$z\mapsto (z,\ell(z))$ as a discrete subset of the 
plane.
Let
$z =\frac{k}{2^n}\in Z_n$ with $n \ge 1$, then we define
the right and the left neighbors
$l(z) =\frac{k-1}{2^n}$ and
$r(z) =\frac{k+1}{2^n}$.
We see that $\ell(l(z)),\ell(r(z)) < n$ and clearly
$l(z)<z<r(z)$, where we take the usual ordering induced by
the reals.
Given $z \in Z$ with $\ell(z) \ge 1$ we consider the
{\em right path}
$z,r(z),r^2(z),\ldots$ and the left path
$z,l(z),l^2(z),\ldots$.
Note that after a finite number of steps the right
path always ends at $1$ and the left path always
ends at $0$.

We use the following facts:

Fact 1: 
Consider for an arbitry 
$z\in Z$ the levels of the vertices on the right and on the left path, i.e.
$\ell(l(z),\ell(l^2(z)),\ldots$ and
$\ell(r(z),\ell(r^2(z)),\ldots$. 
Then all intermediate levels
$n$ with $0<n <\ell(z)$ occur exactly once (either on the right or on the
left path).
E.g. consider
$11/64$ which is of level $6$.
The left path is
$11/64,5/32,1/8,0$ (containig the intermediate levels $5$ and $3$),
the right path is
$11/64,3/16,1/4,1/2,1$ (containing the remaining intermediate
levels $4,2$ and $1$).
This fact can be verified by looking to the dyadic expansion
of $z$, e.g. $11/64 =0.001011$. Note that the dyadic expression 
of $l(z)$ is
obtained from the one of $z$ by removing the last $1$ in this expression,
i.e. $l(0.001011)=0.00101$. The dyadic expression of $r(z)$ is obtained
by removing the last consequetive sequence of $1$'s and putting a $1$ 
instead of the $0$ in the last entry before the sequence, e.g
$r(0.001011)=0.0011$. Therefore the levels of the left path (resp. of 
the left path)
correspond to the
places with a $1$ (resp. with a $0$) in the dyadic expansion.

Fact 2: Let $l^k(z)$ be a point on the left path and
$\ell(l^k(z))\ge 1$. Let $m$ be the integer, such that
$r^m(z)$ is the first point on the right path with
$\ell(r^m(z)) < \ell(l^k(z))$,
then
$r(l^k(z))=r^m(z)$.
A corresponding statements holds for points on the right path.
This fact can also be verified by looking to the dyadic expension.

We consider the graph whose vertex set is $Z$, and the edges
are given by the pairs
$\{0,1\} ,\{z,r(z)\}, \{z,l(z)\}$, where the
$z \in Z$ are points with level $\ge 1$.
One can visualize this graph nicely, if we use
the realization of $Z$ in the plane described above.
In this picture we can see the edges as line intervalls
and the graph is planar.
In this picture the left path from a point
$z$ with $\ell(z) \ge 1$ can be viewed as the graph of
a piecewise linear function defined on the
interval $[0,z]$ (here $z\in [0,1]$) and the the right path
as the graph of a piecewise linear function on $[z,1]$.
The union of these two paths form a "tent" in this picture.

Fact 3: Below this tent there lies not point of $Z$.

%%%%%%%%%%%%%%%%%%%%%%%%%%%%%%%%%%%%%%%%%%%%%%%%%%%%%%%%%%%%%
\begin{figure}[htbp]
\centering
\includegraphics[width=1.0\columnwidth]{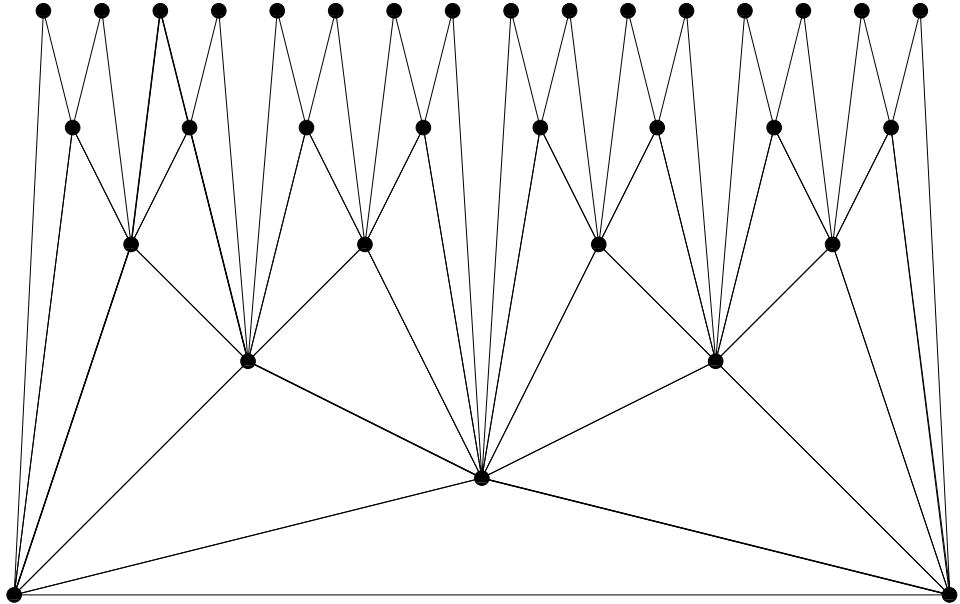}
\caption{graph with a tent}\label{tent}
\end{figure}
%%%%

To every edge in this graph we associate a length.
To the edge $\{0,1\}$ we associate the length $1$,
and to an edge of the type
$\{z,l(z)\}$ and 
$\{z,r(z)\}$ we associate the length
$a^{\ell(z)}$.
Now we define the quasimetric $\rho$.
First set $\rho(0,1)=1$.
Let $z,z' \in Z$ be points such that
$z, z'$ is not the pair $0,1$. Let us assume $z < z'$.
Then we consider the right path
$z,r(z),r^2(z),\ldots,1$ starting from $z$, and
the left path $z',l(z'),l^2(z'),\ldots,0$
starting from $z'$.
Then the properties from above imply that 
these two paths intersect at a unique
point $r^k(z)=l^s(z')$.
Then we obtain a V-shaped path
$z,r(z),\ldots,r^k(z)=l^s(z'),\ldots,l(z'),z'$
formed by edges from our graph from
$z$ to $z'$.
We define
$\rho(z,z')$ to be the sum of the lengths of
the edges of this path.

The main point is now to show that
$\rho$ is a quasi-metric space.
Before we prove this, we show that
the chain approach does not give a metric.
Therefore consider for any integer $n$ the
chain
$0,1/2^n,2/2^n,3/2^n,\ldots 2^n/2^n=1$.
by our definition
$\rho(i/2^n,(i+1)/2^n) = a^n$.
Thus the length of the chain is
$2^n a^n$ which converges to $0$ since 
$a < 1/2$.

It remains to show that
$\rho$ is a quasi-metric.

Fact 4: If $z \in Z_n$, then
$\rho(z,0) + \rho(z,1) = \tau_n$ where
$$\tau_n = a+a^2+\cdots +a^{n-1}+2a^n.$$

To obtain this fact consider the tent
formed by the left path from $0$ to $z$ and
the right path from $z$ to $1$ and consider the levels
of the points on this path.
By Fact 1 all intermediate levels occur exactly once.
Thus the formula comes immediately from the definition of $\rho$.
Note that the $2a^n$ comes from the two edges starting at the 
top point $z$ of the tent.

Note that in the "limit case" $a=1/2$ we have
$\tau_n=1$ for all $n$. For $a< 1/2$ we easily compute
$$2a=\tau_1>\tau_2 >\cdots > \tau_{\infty}=\lim \tau_n=\frac{a}{1-a}.$$

Consider now the following special triangle
$z_0,z_1,z_2$, with the properties,
that:

$z_1$ lies on the left path starting from $z_0$

$z_2$ lies on the right path starting from $z_0$

$z_2$ lies on the right path starting from $z_1$

These conditions imply that
$z_1\le z_0\le z_2$ and
$\ell(z_0)\ge \ell(z_1) \ge \ell(z_2)$.

Let $n= \ell(z_0)$ and $m = \ell(z_1)$.
Then Fact 4 applied to the tents
$0,z_0,1$ and $0,z_1,1$ implies that
$$\rho(z_1,z_o)+\rho(z_0,z_2)-\rho(z_1,z_2) = \tau_n-\tau_m < 0,$$
hence
$\{z_1,z_2\}$ is the longest side of that triangle.

We obtain from the above inequalities in particular that

$$  \rho(z_1,z_2)
- (\rho(z_1,z_o)+\rho(z_0,z_2))  \le \tau_m-\tau_{\infty}
=a^m \frac{1-2a}{1-a},$$
where the last equality is an easy computation.
Since
$\rho(z_1,z_2) \ge a^m$, we obtain
$$\frac{(\rho(z_1,z_o)+\rho(z_0,z_2))}{\rho(z_1,z_2)} \ge 1 - \frac{1-2a}{1-a},$$
and hence
$$\rho(z_1,z_2) \le (1+\ep_a) (\rho(z_1,z_o)+\rho(z_0,z_2)),$$
where $\ep_a \to 0$ for $a \to 1/2$. Actually $1+\ep_a =(1-a)/a$.

We consider now an arbitrary (nondegenerate) triangle.
We number the vertices such that
$z_1<z_0<z_2$.

Consider the V-shaped path from $z_1$ to
$z_2$, let
$\tilde{z}=r^k(z_1)=l^s(z_2)$ be the "lowest" point
on this path. By symmetry of the whole argument we assume
without loss of generality that
$z_0 \le \tilde{z}$.
Now (using Fact 3) we see that the left path staring at $z_0$
will intersect the right path starting in $z_1$.
Let $z'_1$ be the intersection point.
Let $z'_2$ be the first point, where the right path starting
at $z_0$ coincides with the right path starting at $z_1$.
Fact 3 implies that $z'_2\le \tilde{z}$.
Note that now the triangle
$z'_1,z_0,z'_2$ is a special triangle as discussed above.
Further note that
$$\rho(z_1,z_0) = \rho(z_1,z'_1) + \rho(z'_1,z_0),$$
$$\rho(z_0,z_2) = \rho(z_0,z'_2) + \rho(z'_2,z_2),$$
$$\rho(z_1,z_2) = \rho(z_1,z'_1) + \rho(z'_1,z'_2)+ \rho(z'_2,z_2).$$
Therefore we see as above
$$  \rho(z_1,z_2)
\ge (\rho(z_1,z_o)+\rho(z_0,z_2)).$$  

We compute

$$\frac{(\rho(z_1,z_o)+\rho(z_0,z_2))}{\rho(z_1,z_2)} \ge
\frac{(\rho(z'_1,z_o)+\rho(z_0,z'_2))}{\rho(z'_1,z'_2)}
\ge 1 - \frac{1-2a}{1-a},$$
where the last inequality is from the special case.
Thus also in this case we obtain
$$\rho(z_1,z_2) \le (1+\ep_a) (\rho(z_1,z_o)+\rho(z_0,z_2)).$$

%%%%%%%%%%%%%%%%%%%%%%%%%%%%%%%%%%%%%%%%%%%%%%%%%%%%%%%%%%%%%%%%%%%%%%%%%%

\bigskip

\end{document}